\documentclass[12pt]{amsart}
\usepackage{amscd}
\usepackage{verbatim}
\usepackage{amssymb, amsmath, amsthm, amscd,ifthen}
\usepackage[dvips]{graphics}
\usepackage[cp866]{inputenc}
\usepackage{graphicx}
\usepackage{epsfig}

\usepackage{amsmath,amssymb,amscd,}
\usepackage[all]{xy}

\usepackage{amssymb}

\newtheorem{thm}{Theorem}[section]

\newtheorem{ex}[thm]{Example}

\newtheorem{lemma}[thm]{Lemma}

\theoremstyle{definition}

\newtheorem{defn}[thm]{Definition}

\newtheorem{alg}{Algorithm}

\begin{document}

\title{Motion planning  and control of a planar polygonal linkage}

\thanks{ The authors acknowledge useful discussions
with George Khimshiashvili and hospitality of CIRM, Luminy,  where
this research was completed in  framework of ''Research in Pairs''
program in January of 2014. The first author was also supported by
RFBR grant 15-01-02021.  We are also grateful to Troy Baker for a number of  English corrections.}

\author{Gaiane Panina}
\address{Institute for Informatics and Automation, St. Petersburg, Russia,
Saint-Petersburg State University, St. Petersburg, Russia.}

\author{Dirk Siersma}
\address{University of Utrecht, Mathematisch Instituut,
Utrecht, The Netherlands.}

\begin{abstract}
For a polygonal linkage, we produce a fast navigation algorithm on
its configuration space. The basic idea is to approximate the configuration space  by
the vertex-edge graph of its cell decomposition
 discovered by the first author. The algorithm has three
aspects: (1) the number of navigation steps does not exceed 15
(independent of the linkage),
 (2) each step is a disguised flex of a quadrilateral from one triangular configuration to another, which is a well understood type of flex,
 and (3) each step can be performed explicitly by adding some extra bars and obtaining a mechanism with one degree of freedom.
\end{abstract}

\maketitle

\section{Introduction}

In the paper we work with a {polygonal linkage} (equivalently, with
a {flexible  polygon}), that is, with a collection of rigid bars
connected consecutively in a closed chain. We allow any number of
edges and any lengths  assignment,  (under a necessary assumption
of the triangle inequality, which guarantees the closing
possibility).  The flexible polygon lives in the plane and admits
different shapes, with allowed self-intersections. Taken together,
all the shapes form the moduli space of the linkage. In the paper,
we produce a motion planning algorithm (equivalently, a navigation
algorithm) which explicitly reconfigures one shape to another via
some continuous motion.  In the language of the moduli space this
means that we present a path leading from one prescribed point to
another. We not only indicate the path, but also present a way of
forcing the linkage to follow the path.

Although the problem does not seem very complicated (since the more
edges we have, the more degrees of freedom  we have), the navigation
is not an easy issue because of the (possible) topological
complexity of the moduli  space.

There exists   (see \cite{LenWh})  an $O(n)$  algorithm, where each
step is a \textit{ line-tracking} motion. That is, during each step,
 the entire polygon except for some pentagonal subchain is frozen,
 only the subchain flexes in such a way that
 one of its vertices   moves along a
straight line.

Our reconfiguring algorithm is  based on a stratification of the
moduli space into a cell complex, introduced in \cite{Pan}. More
precisely, we treat the one-skeleton of the complex as an
appropriate approximation of the moduli space. In other words, we
have an embedded graph, and we mostly navigate along the graph. The
navigation goes as follows:  from a given configuration of the
linkage, we first reach an appropriate vertex of the graph, then
navigate along the graph until we are close to the target
configuration, and next, we pass to the target configuration. There
are three important aspects about the algorithm:
\begin{enumerate}
  \item The number of steps (i.e., the number of edges of  the connecting path) never exceeds 15. That is, we have
  a finite time algorithm   (rather than $n$ or even $log \ n$ complexity).
  \item However, finding each of the 15 designated configurations
requires a linear time complexity algorithm.
  \item Each of the steps (that is, going along an edge of the graph)
   is a disguised  flex of a quadrilateral polygonal linkage, which is both simple and well-understood.
  \item Each of the steps can be performed  explicitly by adding some extra bars and obtaining a mechanism with one degree of freedom, see Section \ref{SecNavigationOnModuli}.
\end{enumerate}

The paper is organized as follows. Section \ref{SecPrelim}  gives
precise definitions and explains the cell structure on the
configuration space. We  also present introductory examples and give
a formula for the number of vertices of the vertex-edge graph of the complex $\Gamma$. Section \ref{SecNavigationOnGraph} explains the
navigation on the graph. We
show that a vertex-to-vertex  navigation requires at most 15 steps.

Our next goal is to realize the prescribed motions. We work under
assumption that we have a full control of convex configurations, that is, we know how to reconfigure one convex configuration to another.
There are different approaches how to  do this: by using Coulomb
potential, as in  \cite{KHPS},  by mechanically controlling the
angles, in the way  described in \cite{ADEHOST}, or in some other
way, not to be discussed in the paper. However we stress that for
navigating over edges of the graph, it suffices to control just
quadrilaterals, which is a much easier task and  well
understood in all respects.

Navigation from an arbitrary point of the moduli space to a vertex
requires one more step: we need to connect the starting point to the
graph. Two different ways of doing this are described in Section \ref{SecNavigationOnModuli}.

The results presented in this paper arose as a natural continuation
of the research on Morse functions on moduli spaces of polygonal
linkages started in \cite{khim2010}, \cite{khipan}, and \cite{KHPS}.
Several approaches to navigation and control of mechanical linkages
have previously been discussed, in particular, in \cite{hausmann},
\cite{hausmannrodriguez},   \cite{khim2010} and \cite{KHPS}.

\newpage
\section{ Moduli spaces of planar polygonal linkages}\label{SecPrelim}
We start by a short review of some results on polygonal linkages and their moduli spaces.

A \textit{polygonal  $n$-linkage} is a sequence of positive numbers
$L=(l_1,\dots ,l_n)$. It should be interpreted as a collection of
rigid bars of lengths $l_i$ joined consecutively in a chain by
revolving joints.  In the literature, it is sometimes called
\textit{ a closed chain}.

We assume that  the closing condition holds: the length of each bar
is less than the total length of the rest.

 \textit{A configuration} of $L$ in the
Euclidean plane $ \mathbb{R}^2$  is a sequence of points
$R=(p_1,\dots,p_{n}), \ p_i \in \mathbb{R}^2$ with
$l_i=|p_i,p_{i+1}|$, and $l_n=|p_n,p_{1}|$.

\begin{defn}
The set $M(L)$ of all  configurations modulo orientation preserving
isometries of $\mathbb{R}^2$ is \textit{the moduli space}, or
\textit{the configuration space, }of the  linkage $L$.
\end{defn}
\begin{defn}\label{defn2}
Equivalently, $M(L)$ can be defined  as the quotient space
$$M(L)=\{(u_1,...,u_n) \in (S^1)^n : \sum_{i=1}^n l_iu_i=0\}/SO(3).$$
\end{defn}

Now let us treat the lengths $l_i$ as variables.  Each $n$-tuple  of positive numbers gives us  a polygonal linkage  which
comes together with its configuration space. The hyperplanes
$$\sum_{i \in J} l_i = \sum_{i\notin J} l_i,$$
where $J$ ranges over all subsets of $[n]$, are called
\textit{walls.} (Here and in the sequel, $[n]$ denotes the set
$\{1,...,n\}$.) The walls decompose $\mathbb{R}_{>0}^n$ into a
collection of {\it chambers}.

 Here is a (far from complete) summary of facts
about $M(L)$:

 \begin{itemize}
\item If no configuration of $L$ fits a straight line, or, equivalently, $L$ does not lie on a wall, then $M(L)$ is a smooth $(n-3)$-dimensional manifold.
In this case, the linkage is called \textit{generic}
\cite{KapovichMillson}. Throughout the paper we consider only
generic linkages.
 \item The topological type of the moduli space $M(L)$ depends only on the chamber containing $L$ \cite{KapovichMillson}.

 \item   $M(L)$ admits a structure of a regular cell complex \cite{Pan}.
  The combinatorics is very much related (but not identical) to the combinatorics of the permutohedron.
  The construction is explained  later in this section.
  \item Definition \ref{defn2} implies that the moduli space $M(L)$, as well as the cell structure,
  does
  not depend on the permutation of the edgelengths $l_i$. More
   precisely, for any permutation $\sigma \in S_n$, there exists
   a canonical isomorphism between $M(L)$ and $M(\sigma(L))$, which
  maintains the cell structures.

 \end{itemize}

\begin{defn}
\label{defn:short} A set $I\subset [n]$ is called {\it long}, if
$$|I|=\sum_{i \in I} l_i > \sum_{i\notin I} l_i.$$
Equivalently, for a long set $I$,
$$|I|>\frac{|L|}{2}.$$

Otherwise, $I$ is called {\it short}.

Note that because of the genericity assumption, we never have
$|I|=\frac{|L|}{2}.$
\end{defn}
%\begin{itemize}
% \item  Homology groups of $M(L)$ are free
% abelian groups. For a generic length vector $L$,
% the rank of the homology group $H_k(M(L))$  equals  $a_k+a_{n-3-k}$, where $a_i$ is the number of short
% sets of cardinality $i+1$ containing the longest edge (see \cite{faS}).
%\end{itemize}

We stress that the manifold  $M(L)$  (considered either as a
topological manifold, or as a cell complex) is uniquely defined by
the collection of short subsets of $[n]$.

\bigskip

\begin{ex}\label{exBow} Assume that for an $n$-linkage $L$, we have the following:
 $$\forall \ i=1,2,...,(n-1), \hbox{ the set } \{n,i\} \hbox{ is long.}$$
 Loosely speaking, we have ''one long edge''. Such a linkage we call an
 $n$-bow;
its moduli space is a sphere (see \cite{faS}).
\end{ex}

\begin{defn}
A partition of the set $\{l_1,\dots ,l_n\}$ is called
\textit{admissible } if all the parts are non-empty and short.

Instead of partitions of $\{l_1,\dots ,l_n\}$ we shall speak of
partitions of the set $[n]$, keeping in mind the lengths
$l_i$.
\end{defn}

\subsection*{Description of the cell complex}
Now we sketch the cell complex structure on  the moduli space
$M(L)$, referring the reader to \cite{Pan} for even more details and all the
proofs.

\bigskip

The cell structure comes from the following \textit{labeling} of
configurations. Given a configuration $P$ of $L$ \textbf{without
parallel edges}, there exists a unique convex polygon $Conv(P)$
which we call \textit{the convexification of} $P$ such that

\begin{enumerate}
    \item The edges of $P$ are in one-to-one correspondence with the edges of
    $Conv(P)$ . The bijection preserves the directions of the
    vectors.

    \item The  edges of
    $Conv(P)$  are oriented in the  counterclockwise direction with respect to
    $Conv(P)$ .
\end{enumerate}

\begin{figure}[h]
\centering
\includegraphics[width=10 cm]{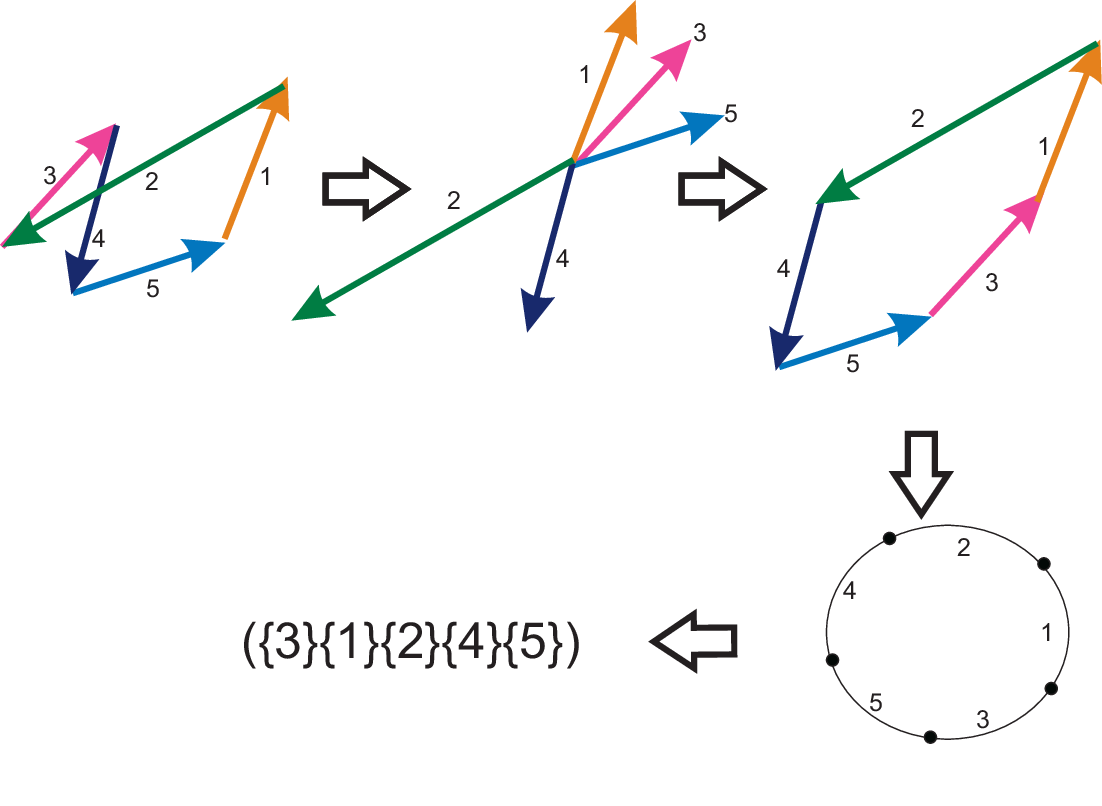}
\caption{Labeling of a polygon with no parallel
edges}\label{rearrangment1}
\end{figure}

In other words, the edges of  $Conv(P)$  are the edges of $P$
ordered by the slope (see Fig. \ref{rearrangment1}). Obviously,
$Conv(P) \in M(\lambda L)$ for some permutation  $\lambda
 \in S_n$. The permutation is
defined up to some power of the cyclic permutation
$(2,3,4,...,n,1)$, and hence can be considered as a cyclic ordering
on the set $[n]$.

\begin{figure}[h]
\centering
\includegraphics[width=8 cm]{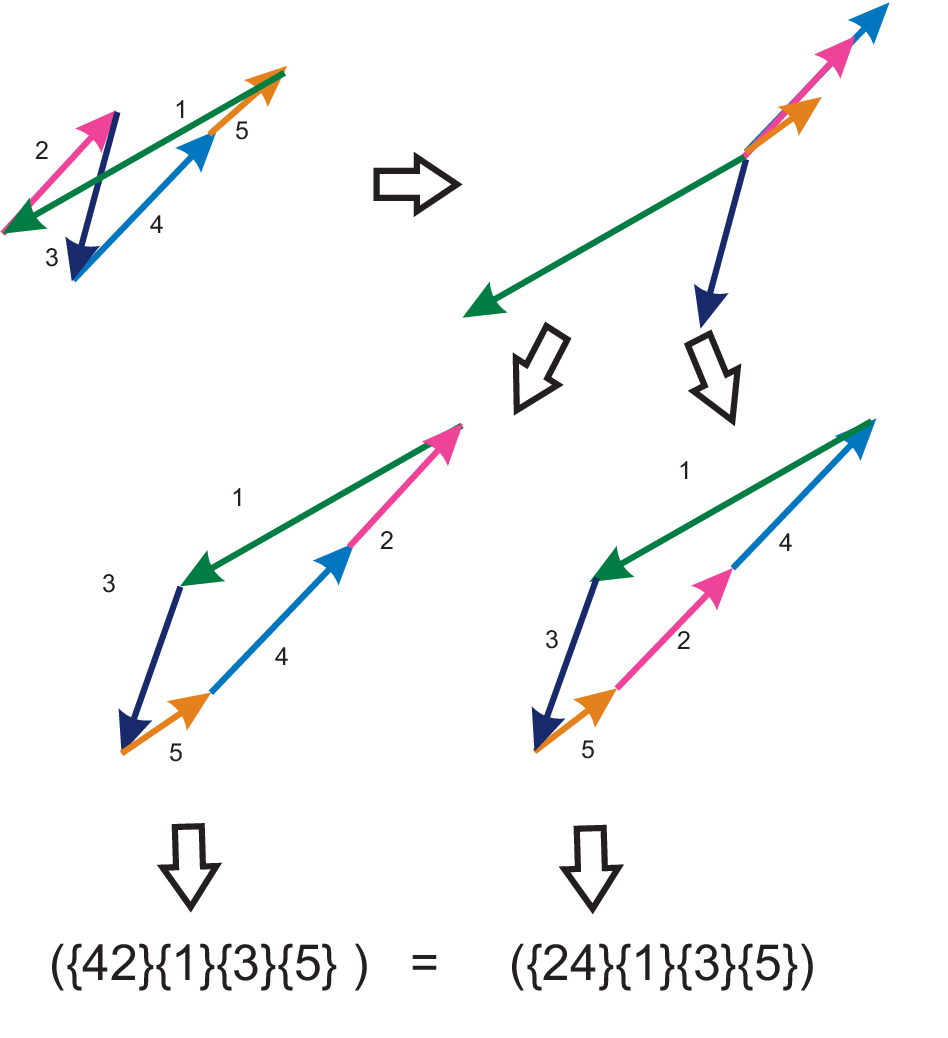}
\caption{Labeling of a polygon with parallel
edges}\label{rearrangment2}
\end{figure}

%Conversely, each convex polygon from $ M(\lambda L)$  is the
%image of some element of $M(L)$ under the above rearranging map.

Our construction assigns to $P$ the \textit{label} $\lambda$,
considered as a cyclic ordering  on the set $[n]$. Equivalently, a
label of a configuration without parallel edges is a cyclically
ordered partition of the set $[n]$ into $n$ singletons. Figure
\ref{rearrangment1} gives an example.
\bigskip

If $P$ \textbf{has parallel edges}, a permutation which makes $P$
convex is not unique, since there is no ordering on the set of
parallel edges. The label assigned to $P$ is a cyclically ordered
partition of the set $[n]$, where parallel edges  belong to the same subset in the partition. Fig. \ref{rearrangment2} gives an
example.

An obvious observation is that all labels are admissible partitions.
\bigskip

{\bf A remark on notation.} We write a cyclically ordered partition
as a (linearly ordered) string of sets, keeping in mind that the
string is closed.

 We stress once again that
there is no ordering inside a set, that is, two labels are equal when they differ only by
permutations of the elements inside the sets. For instance,
$$(\{1\} \{3 \} \{4  2 56\})\neq(\{3 \}\{1\}  \{4  2 56\})= ( \{3 \}\{1\}\{ 24 56\}).$$

\begin{defn}
Two points from $M(L)$ (that is, two configurations) are
\textit{equivalent} if they have  the same label. Equivalence
classes  of $M(L)$   are the \textit{open cells}. The closure of an
open cell in  $M(L)$ is called a \textit{closed cell}. For a cell
$C$, either closed or open, its label $\lambda (C)$ is defined as
the label of any interior point.

\end{defn}

With this labeling, the following theorem is valid.
\begin{thm}\label{MainThm}\cite{Pan}
The above described collection of open cells yields  a structure of
a regular CW-complex $\mathcal{K}(L)$ on the $(n-3)$-dimensional
manifold $M(L)$. Its complete combinatorial description reads as
follows:
\begin{enumerate}
    \item  $k$-cells of the complex $\mathcal{K}(L)$ are labeled by cyclically ordered admissible partitions of
 the set  $[n]$  into $(k+3)$ non-empty
parts.
\item
In particular, the facets of the complex (that is, the cells of
maximal dimension)
    are labeled by cyclic orderings of the set $[n]$.
    \item A closed cell $C$ belongs to the boundary of some other closed cell
    $C'$  iff  the partition  $\lambda(C')$ is finer than  $\lambda(C)$.
    \item The complex is regular, which means that each $k$-cell is attached to the $(k-1)$-skeleton
    by an injective mapping. All closed cells are ball-homeomorphic.
    \qed

\end{enumerate}
\end{thm}

\begin{ex}\label{ExTwoVertices}
The vertex labeled by $$(\{1,2,5,6\},\{3,4\},\{7,8\})$$  and the
vertex labeled by $$(\{1,2\},\{3,4,5,6\},\{7,8\})$$ are connected by
an edge  labeled by   $$(\{1,2\},\{5,6\},\{3,4\},\{7,8\}).$$
\end{ex}
\bigskip

\begin{figure}[h]
\centering
\includegraphics[width=6 cm]{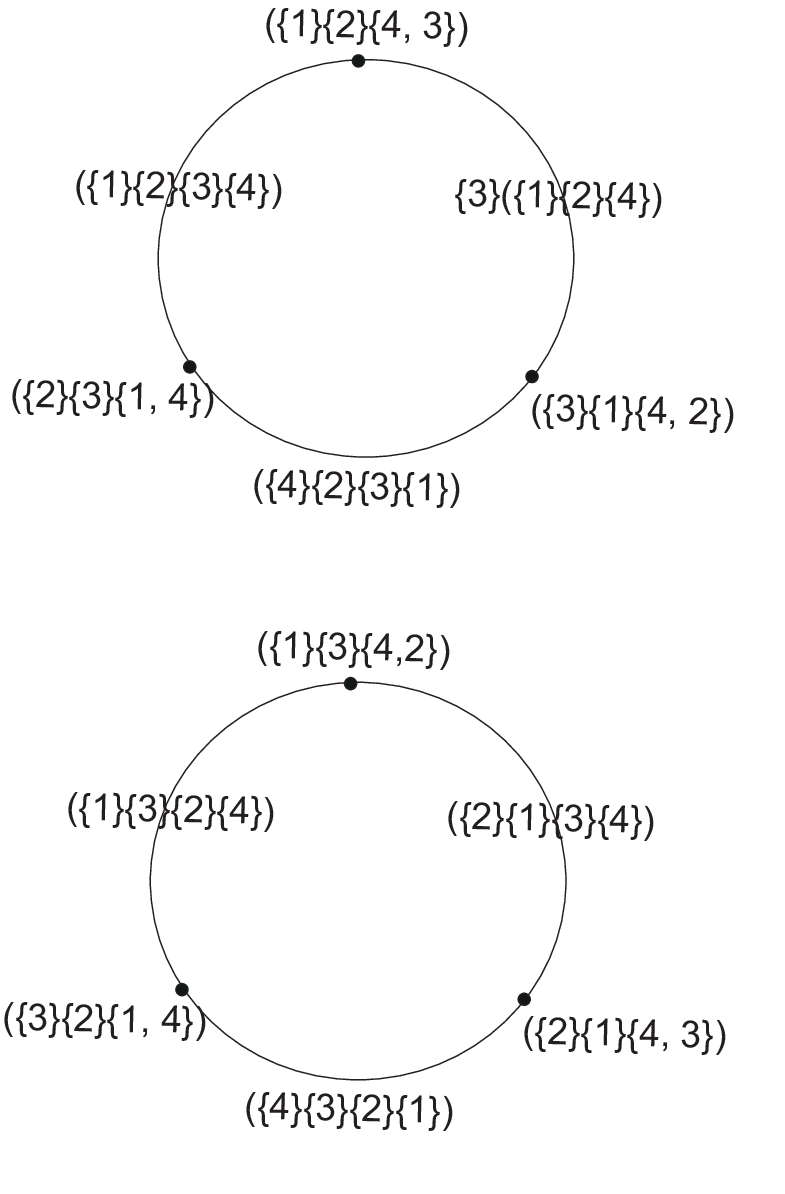}
\caption{$\mathcal{K}(L)$ for the 4-gonal linkage
$(1,1,1,1/2)$}\label{4gon1}
\end{figure}

\begin{ex}  Let $n=4; \ \ l_1=l_2=l_3=1,\ l_4=1/2.$
The moduli space $M(L)$ is known to be a disjoint union of two
circles (see \cite{faS}). The cell complex $\mathcal{K}(L)$ is
depicted in Fig. \ref{4gon1}.
\end{ex}

\begin{ex} \label{ExHexagon}  Assume that for a $4$-linkage $L$,
the sets $\{2,3\}$, $\{4,3\}$, and $\{2,4\}$ are short. An example
of such a length assignment is $$ l_4=l_2=l_3=1,\ l_1=2.5.$$

The moduli space $M(L)$ is  homeomorphic to a  circle.
 The cell complex is combinatorially a hexagon, that is, there are six vertices connected by six edges into a circle.
 The (cyclic) order of the labels of the six vertices is:
 $$(\{1\}\{2,3\}\{4\})$$
 $$(\{1\}\{2\}\{4,3\})$$
 $$(\{1\}\{2,4\}\{3\})$$
 $$(\{1\}\{4\}\{2,3\})$$
 $$(\{1\}\{4,3\}\{2\})$$
 $$(\{1\}\{3\}\{4,2\}).$$

\end{ex}

\begin{ex}For the equilateral pentagonal linkage $L=(1,1,1,1,1)$,
the complex $\mathcal{K}(L)$ has $30$ vertices, $60$ edges, and $24$
pentagonal $2$-cells. Each vertex is incident to exactly four edges.
\end{ex}

\bigskip

In the paper we shall make use of the \textit{vertex-edge graph}
$\Gamma$  of the cell complex, that is, we take into account only
zero- and one-dimensional cells.  We treat it as a (combinatorial)
graph, also keeping in mind its embedding in the $M(L)$. This allows
us to view the graph $\Gamma$  as a discrete approximation of the
moduli space.

\subsection*{Combinatorics of the vertex-edge graph}  Assume that a linkage $L$ is fixed.
As a particular case of Theorem \ref{MainThm}, vertices of $\Gamma$
are labeled by   cyclically ordered partitions of $[n]$ into three
non-empty short sets, and the edges are labeled by cyclically
ordered partition of $[n]$ into four non-empty short sets.  Two
vertices labeled by $\lambda$ and $\lambda'$ are joined by an edge
whenever the label $\lambda$ can be obtained from  $\lambda'$  by
shifting some amount of entries from one set to another, as in
Example \ref{ExTwoVertices}.

\subsection*{Embedding of the vertex-edge graph}

Now we describe a way of retrieving the vertex-edge graph $\Gamma$ from the labels.

\begin{figure}[h]
\centering
\includegraphics[width=10 cm]{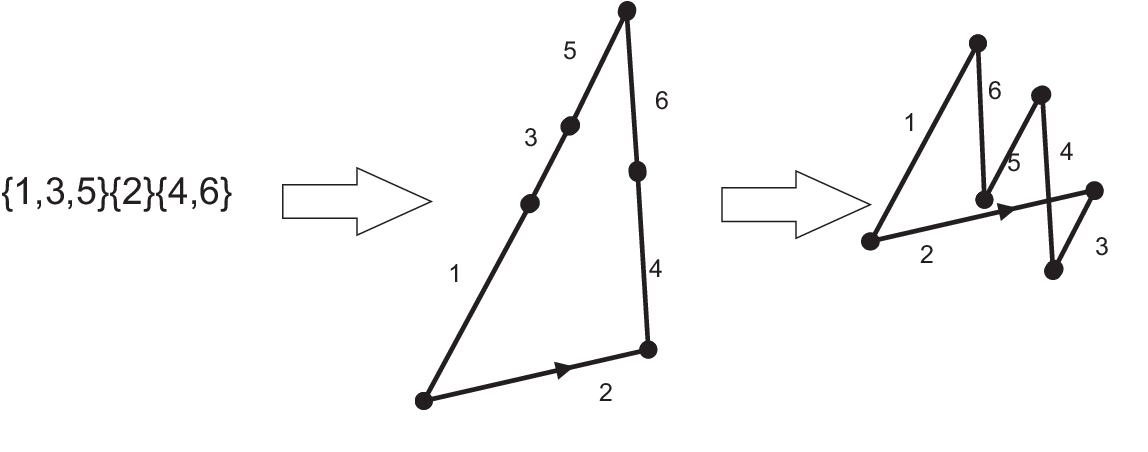}
\caption{Retrieving  a vertex from its label. Each vertex is a
disguised triangle}\label{FigVertex}
\end{figure}

\begin{alg}\textbf{ Retrieving a vertex (viewed as a polygon) by its
label.}

Given a label $\lambda=(I,J,K)$, it  corresponds to a unique point
$P \in M(L)$, that is, to some configuration of $L$. The polygon $P$
can be constructed as follows.
\begin{enumerate}
  \item Take a positively oriented triangle with edgelengths $$\sum_Il_i,\ \sum_Jl_i, \hbox{ and }\sum_Kl_i.$$
  \item We assume that each edge  is composed of segments $l_i$. For instance, we
   decompose the first edge into segments of lengths $\{l_i\}_{i\in I}.$ Their order does not matter.
  \item Now take all the edges apart, keeping their directions, and rearrange them according to the ordering $1,2,...,n$.
  We get a closed polygonal chain $P$, see Figure \ref{FigVertex}.
\end{enumerate}

\end{alg}

Taken together, all labels give all configurations that have exactly
3 directions of the edges. They form the vertex set of the embedded
graph $\Gamma(L)$.

Now let us explain  how the edges are embedded.

\begin{figure}[h]
\centering
\includegraphics[width=14 cm]{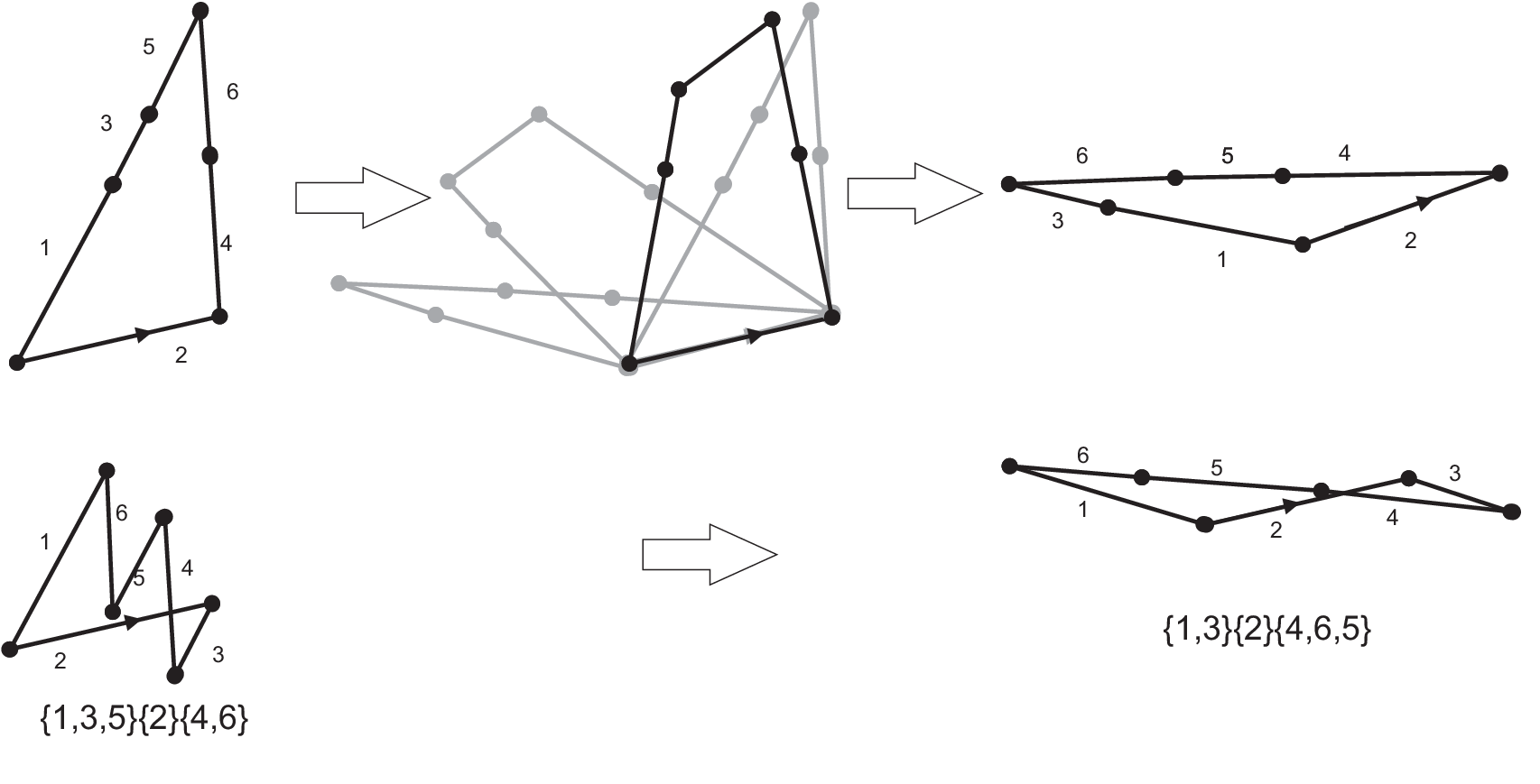}
\caption{Retrieving  an edge from the label
$(\{5\}\{1,3\}\{2\}\{4,6\})$. Each edge is a disguised flex of a
convex quadrilateral from one convex triangular configuration to
another }\label{FigEdge}
\end{figure}

\begin{alg} \textbf{Retrieving an edge by its
label}

Given a label $\lambda=(I,J,K, N)$, it labels an embedded edge of
$\Gamma(L)$, that is, a one-parametric family of configurations.
They are retrieved as follows.
\begin{enumerate}
  \item Take a positively oriented  convex quadrilateral  with consecutive edgelengths $$\sum_Il_i,\ \sum_Jl_i\ \sum_Kl_i, \hbox{ and } \sum_Nl_i.$$
  The set of such quadrilaterals forms a path in $M(L)$  going from one triangle to another, see Figure \ref{FigEdge}.
  \item As in Algorithm 1, decompose each edge  into (short) edges $l_i$.
  \item Exactly as in algorithm 1, take all the edges apart, keeping their directions, and rearrange them according to the ordering $1,2,...,n$. This gives
   a one-parametric family of closed polygonal chains.
\end{enumerate}

\end{alg}

\bigskip

In other words, each embedded edge of the graph $\Gamma(L)$
corresponds to a flex of some quadrilateral, which connects two
triangular configurations. Such a flex can be performed by
compressing one of the diagonals.

\begin{lemma} \begin{enumerate}
                \item The number of vertices of $\Gamma(L)$ equals
$$\sum_{k=1}^n N_k2^{n-k} -2\cdot3^{n-1}+2^n,$$
where $N_k$ is the number of short $k$-sets.
                \item For the $n$-bow, $\Gamma(L)$ has the minimal possible
                number of vertices among all $n$-linkages. In this case, it equals $2^{n-1}-2$.
              \end{enumerate}

\end{lemma}
Proof. We start with  the bow linkage (see Example \ref{exBow}), and
then change $L\in \mathbb{R}^n$ continuously. From the chamber that
corresponds to the bow, we can reach any other chamber by crossing
walls. As follows  from the construction of the cell complex, the
graph $\Gamma$ changes only when $L$ crosses a wall. ''Crossing a
wall'' means that exactly one short $k$-set $I$ becomes long, and
exactly one long $(n-k)$-set $\overline{I}$ becomes short, where the
number $k$ depends on the wall.
 Observe that any proper subset of $\overline{I}$
is short both before and after crossing a wall. A vertex is eliminated
whenever its label contains $I$. A vertex arises whenever its label
contains  $\overline{I}$. This means that after crossing a wall, the
number of vertices changes by adding
$(2^k-2)-(2^{n-k}-2)=2^k-2^{n-k}$.

Therefore,
$$|Vert(\Gamma)|=\sum_{k=1}^{n}N_k2^{n-k}+X_n,$$
where $X_n$ depends solely on $n$. Our second observation is that
for a bow  linkage, labels of the vertices are as follows:
$$( I, \ [n-1] \setminus {I}, \{n\}),$$
where $I$ is  any proper non-empty  subset of $ [n-1]$. Therefore,
$$|Vert(\Gamma(n\hbox{-bow}))|=2^{n-1}-2.$$  This is the base case of an
induction.

For the $n$-bow,$$N_k=\left\{
                        \begin{array}{ll}
                          0=C_{n-1}^0-1, & \hbox{if } k=0; \\
                          n=C_{n-1}^1+1, & \hbox{if } k=1; \\
                          0=C_{n-1}^{n-1}-1, & \hbox{if }  k=n-1 \\
                          C_{n-1}^k, & \hbox{ otherwise.}
                        \end{array}
                      \right.
$$

$$2^{n-1}-2=\sum_{k=0}^{n}C_{n-1}^k2^{n-k}+ 2^{n-1}-2^n-2+X_n$$
$$X_n=2^n-2\cdot3^{n-1},$$
which yields the final formula.\qed

\bigskip

As we see, the number of vertices of the graph $\Gamma$ depends exponentially on $n$.
However, the valence of the vertices also depends on $n$ exponentially, so one can expect a
small diameter (in the graph-theoretic sense, that is, the maximal length of the shortest path between two vertices), and a fast navigating algorithm.

\begin{lemma}A vertex of the graph $\Gamma$ labeled by $(I,J,K)$ has exactly
$$2^{i}+2^{j}+2^{k}-6$$ incident edges,
where $i,\ j, $   and $k$ are the numbers of elements in $I,\ J,$
and $K$ respectively.\qed
\end{lemma}

\section{ Motion planning on the graph $\Gamma(L)$}\label{SecNavigationOnGraph}
Here we describe a  finite-step algorithm of motion planning (or,
equivalently, navigation) from an arbitrary vertex of $\Gamma(L)$ to
an arbitrary target  vertex.

A \textit{path}  means a graph-theoretical path, that is, a
consecutive  collection of edges. By its \textit{length}, or
\textit{the number of steps} we mean just the number of edges in the
path.

We start with an example demonstrating that the proposed navigation
works quickly.

\begin{alg}\label{AlgBow} \textbf{(Navigation for the bow linkage)} For   an $n$-bow linkage with $n>3$,
any two vertices of $\Gamma(L)$ are connected by a path whose length is at most 3.
The path is explicitly described as follows.
\begin{enumerate}
  \item We start with a vertex labeled by $$(I,J,\{n\}).$$  We may  assume that the target vertex is labeled by
  $$(\{1,2,...,k\},\{k+1,k+2,...,n-1\},\{n\}).$$
  If this is not the case, we can renumber the edges of the linkage:
  we know that
  in view of Definition \ref{defn2}, renumbering maintains the manifold $M(L)$ and the cell
  structure.
  \item If  $|J|>1$,
  \begin{enumerate}
  \item If $1 \in I$, then go to the step (b).
  If not, shift  the entry $1$  to the set $I$ and go to step (b).
  \item Shift the set $I \setminus \{1\}$ to the set $J$. We arrive at the vertex labeled by  $$(\{1\},\{2,3,...,n-1\},\{n\}).$$
  \item Shift  the set $\{2,3,...,k\}$  to the first set and get the target vertex.\qed
\end{enumerate}
 \item If  $|J|=1$, then $|I|>1$. We act similarly to the
 case (2), exchanging the roles of $I$ and $J$. Namely:
\begin{enumerate}
  \item If $(k+1) \in J$, then go to the step (b).
  If not, shift  the entry $(k+1)$  to the set $J$ and go to step (b).
  \item Shift the set $J \setminus \{k+1\}$ to the set $I$. We arrive at the vertex labeled by  $$(\{1,2,...,k,k+2,...\},\{k+1\},\{n\}).$$
  \item Shift  the set $\{k+2,k+3,...,n-1\}$  to the second set.\qed
\end{enumerate}
\end{enumerate}
\end{alg}

Now we pass from a bow to arbitrary linkages. We first produce an
auxiliary algorithm which\textit{ turns polygons inside out}, that
is,  connects a configuration with its mirror image by a path.
\begin{alg}\label{alg5}\textbf{(Turning a pentagon inside out)} Assume that a $5$-linkage $L$ satisfies the following conditions:
\begin{enumerate}
  \item The set $\{1,2\}$ is long.
  \item The set $\{1,5\}$ is long.
  \item $\forall \ i\neq1, \ j\neq1 ,$ the set $\{i,j\}$ is short.
\end{enumerate}
Then there exists a 4-steps path which turns the configuration
$(\{4,5\}\{1\}\{2,3\})$ inside out. Here it is:
 $$(\{4,5\}\{1\}\{2,3\})$$
 $$(\{4,5,3\}\{1\}\{2\})$$
 $$(\{4,3\}\{1\}\{2,5\})$$
 $$(\{4,3,2\}\{1\}\{5\})$$
 $$(\{3,2\}\{1\}\{4,5\})$$
\end{alg}

\bigskip

\begin{alg}\label{AlgINsideOut} \textbf{(Turning an arbitrary polygon inside out)}
Assume that a linkage $L$ is fixed.
\begin{enumerate}
  \item If the configuration space  $M(L)$  is connected, then from each vertex $V$   labeled by $(I,J,K)$ there are at most
  six
  steps
 to its mirror image $(K, J, I)=(J,I,K)$.
  \item  If the $M(L)$ is disconnected, then the vertex $(I,J,K)$  and  its mirror image $(J,I,K)$ lie in different connected components, and no connecting path exists.
\end{enumerate}
The idea is to imitate   a pentagon which satisfies the three conditions of the Algorithm \ref{alg5} by freezing
 some of the entries in $L$. Here and in the sequel, ''freezing'' means putting the entries into one separate set and after that, manipulating with the set as with a single  entry.

 Assume that $l_1\geq l_k\geq l_m$ are the longest edges of $L$.
 It is known from \cite{KapovichMillson} that $M(L)$ is connected if and only if $$l_k+l_m< \frac{|L|}{2}.$$

 Our starting point is a vertex labeled by $(I,J,K)$, where $1 \in
 J$. We can assume this, since $(I,J,K)=(J,K,I)=(K,I,J).$
 Since we can apply renumbering of the edges, we also can assume that the entries are as
 follows:
 $$(I,J,K)= (\{r+1,r+2,...,p\}\{p+1,p+2,...,n,1,2,...,q\}\{q+1,q+2,...,r\}).$$

\textbf{Steps 1--2: necessary preparations before freezing: pushing
entries to the central set.}

 Maintaining the ordering, we shift to the middle set as many entries from the first and the last  sets as is possible. In more
 detail,
 we first decide what entries should be shifted from $I$, and what entries should be  shifted from $K$. After the decision is done,
  we make two shifts, which means two steps. The choice is not uniquely defined; however,  any choice will
suffice.
 The result we denote by $$(I\setminus I',J\cup I'\cup K',K\setminus
 K'),$$
for which we keep the same notation:

$$(I\setminus I',J\cup I'\cup K',K\setminus
 K')=$$ $$=(\{r+1,r+2,...,p\}\{p+1,p+2,...,n,1,2,...,q\}\{q+1,q+2,...,r\}).$$

%Observe that after shifting,
 %the sets $$\{p,p+1,p+2,...,n,1,2,...,q\} \hbox{ and  }\{p+1,p+2,...,n,1,2,...,q, q+1\} \hbox{ are long. }$$

\bigskip

\textbf{Steps 3--6: freeze the linkage either to a quadrilateral or
to a pentagon, and then turn it inside out.} On this step we work
under assumption that $M(L)$ is connected.
 A necessary observation
is  that now the set $\{p,q+1\}$ is short.
 Indeed, $$l_p+l_{q+1}\leq l_k+l_m < \frac{|L|}{2}.$$  Therefore  the set $$A=\{q+2,q+3,..,r,r+1,...,p-2,p-1\}$$ is not empty.

Now the algorithm splits depending on the number of entries in $A$
which equals $(p-1)-(q+2)+1=p-q-2$.

 \begin{enumerate}
   \item Assume that  $p-q-2=1$, which means that $A$ is a one-element set. Then
 we freeze the  set
 $\{p+1,p+2,...,n,1,2,...,q\}.$
After renumbering
 $$4:= \{p\}, \ 1:=\{p+1,p+2,...,n,1,2,...,q\},$$$$ \ 2:=\{q+1\}, \hbox{ and } 3:=A,$$
 we arrive at a quadrilateral from Example \ref{ExHexagon} which  can be turned inside out in three steps.

   \item   Assume that  $p-q-2>1$. We divide the set $A$ into two non-empty subsets and freeze the two subsets.
   We also freeze the set $\{p+1,p+2,...,n,1,2,...,q\}.$
   That is, for instance,
  we freeze the following three sets of entries:
 $$\{r+1,r+2,...,p-1\},\ \{p+1,p+2,...,n,1,2,...,q\}, \  \hbox{ and } \{q+2,...,r\}.$$
After renumbering
 $$4:=\{r+1,r+2,...,p-1\},\ 5:= \{p\}, \ 1:=\{p+1,p+2,...,n,1,2,...,q\},$$$$ \ 2:=\{q+1\}, \hbox{ and } 3:=\{q+2,...,r\},$$
 we arrive at a pentagon which satisfies the properties from the Algorithm \ref{alg5}, and therefore can be turned inside out in four steps.
 \end{enumerate}
\end{alg}
\textbf{Steps 7--8. Pushing entries from the central set.} We have
arrived at a vertex labeled by $$(K\setminus
 K',J\cup I'\cup K',I\setminus
I').$$ In two steps we get to $(K,J,I)$. These steps are reverse to
steps 1--2.

\bigskip

\begin{alg} \label{PropEitherOr} For any $n$-linkage  $L$, and any  two vertices $V$ and $V'$ of the graph $\Gamma(L)$, the vertex $V$ is connected either to $V'$ or to the
 mirror image of $V'$ by  a path whose  length is at most 7.
The path is explicitly described   as follows. Assume that $l_1$ is
the longest edge, and that the target vertex $V'$ is labeled by
  $$(\{1,2,...,k\},\{k+1,k+2,...,m\},\{m+1,m+2,...,n\}).$$

\begin{enumerate}
  \item In two steps we get  from $V$ to a vertex  labeled by $$(I,\{1\},J)$$  for some $I$ and $J$. This is always possible:
 \begin{enumerate}
   \item Assume that $V$ is labeled by $(A,B,C)$, and $1\in B$. Start shifting the entries from $B\setminus \{1\}$ to the set $C$. We shift  as many  entries as possible, that is, we think of shifting them one by one until $C$ is short, and stop if no other entry can be added to $C$  without making it long. The order in which we treat the entries does not matter. However, we should shift all the entries by one step, that is, first decide what entries are to be shifted, and next, shift them as  one subset.
   \item Shift the rest of  $B\setminus \{1\}$ to the set $A$.

 \end{enumerate}

  \item If one of the sets $I$ or $J$ contains two consecutive entries, we can freeze them. We freeze all possible pairs of consecutive entries and renumber the edges (preserving the ordering), which gives  us a linkage with a smaller
      number of edges.

       In any case we have a vertex labeled  either by
      $$(\{3,5,7,...\}, \{1\},\{2,4,6,8,...\}),$$
      or by the symmetric image
      $$(\{2,4,6,8,...\}, \{1\},\{3,5,7,...\}).$$
      \item Pull $2,3,4,...,s$ and $3,5,7,...,s\pm 1$ to the middle set for the largest $s$ which is possible. (This requires 2 steps more).
       Thus   we arrive either at  the vertex
           $$(\{s+1,s+3,...\}, \{1,2,3,4,...,s\},\{s+2,s+4,...\}),$$
           or at its symmetric image
           $$(\{s+2,s+4,...\}, \{1,2,3,4,...,s\},(\{s+1,s+3,...\}).$$
      So the first entry that we cannot shift  to the middle set is $s+1$.
      \item Shift $s+3,s+5,...$ to the third set. It is possible because \newline $\{1,2,3,4,...,s,s+1\}$ is long.
      We arrive either at
       $$(\{s+1\}, \{1,2,3,...,s\},\{s+2,s+3,....,n\})=$$$$= (\{s+1\}, \{s,s-1,...,2,1\},\{n,n-1,...,s+2,s+3\}),$$
       or at the symmetric image. Now we have either clockwise or counterclockwise ordering on the entries.
      \item There are two steps either to the target, or to the mirror image of the target vertex. \qed
\end{enumerate}
\end{alg}

Combining the two above algorithms, we immediately have the theorem:
\begin{thm} \begin{enumerate}
              \item For any $n$-linkage  $L$  with a connected moduli space,  any  two vertices $V$ and $V$' of the graph $\Gamma(L)$ are connected  by  a path whose  length is at most 15.
              \item For any $n$-linkage  $L$  with a disconnected moduli space, and any  two vertices $V$ and $V$' of the graph $\Gamma(L)$ lying in the same connected component, $V$ and $V$'
              are connected  by  a path whose  length is at most 7.
            \end{enumerate}
The path is constructed explicitly using the above algorithms.
That is we have  the following algorithm:

\bigskip

\textbf{Algorithm A}
\begin{enumerate}
  \item Starting with a vertex $V$, follow  Algorithm \ref{PropEitherOr}. It may bring us  to the target
  point, and then we are done.
  \item If on the first step we get the mirror image of the target point, turn it inside out via Algorithm \ref{AlgINsideOut}. \qed
\end{enumerate}

 \end{thm}

\begin{figure}[h]
\centering
\includegraphics[width=14 cm]{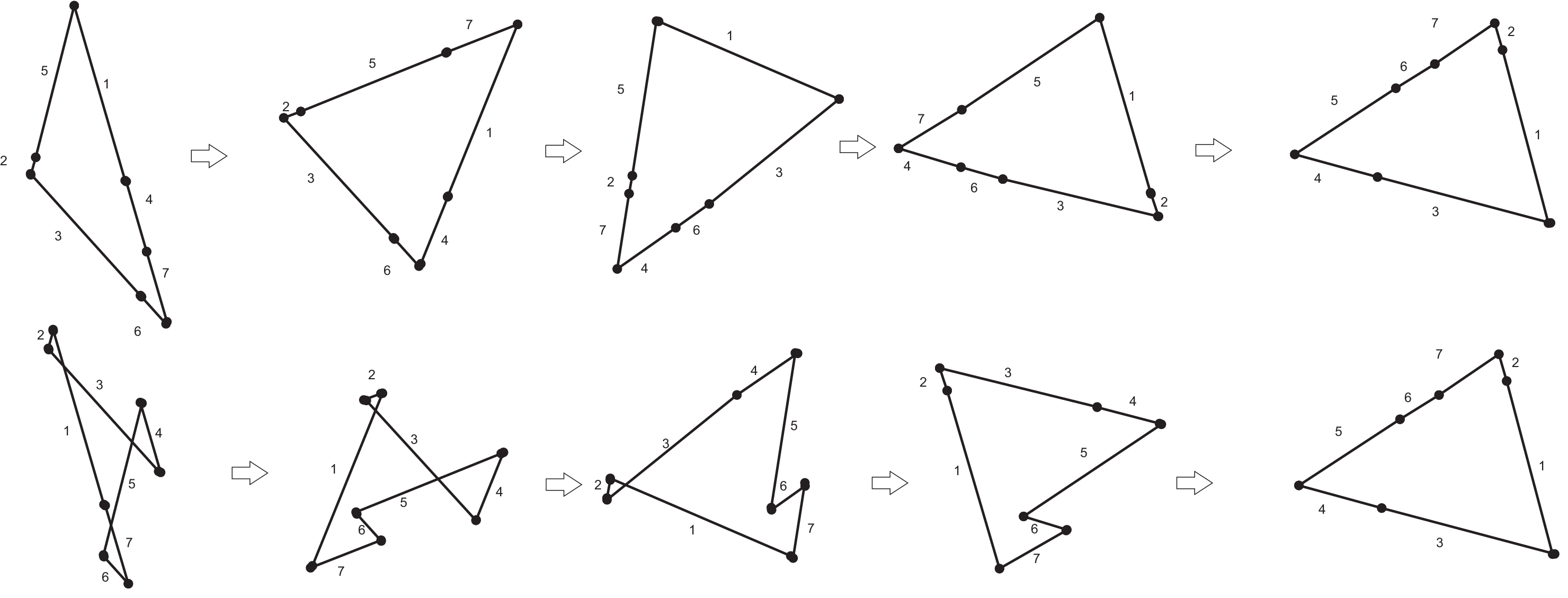}
\caption{The first part of Algorithm A applied to the heptagonal
configuration. The second row depicts the vertices of the path. The
first row depicts the disguised flex. }\label{FigHeptagon}
\end{figure}

Now we exemplify the steps of the algorithm for one particular
heptagonal configuration.
\begin{ex} Assume we have a heptagonal linkage with edge lengths
$$l_1=10, \ l_2=1, \ l_3=9, \ l_4=4, \ l_5=9, \ l_6=2, \  \hbox{ and } l_7= 4.$$

Assume that our starting point is $V_1=(\{3,6\}\{1,4,7\}\{2,5\})$,
and that the target vertex is $V'=(\{5,6,7\}\{1, 2\}\{3,4\})$. Then
Algorithm  A  runs as is described below and as is depicted in
Figure \ref{FigHeptagon}.
\begin{enumerate}
                                     \item The starting point is the vertex of the graph  $$V_1=(\{3,6\}\{1,4,7\}\{2,5\}).$$
                                     According to Algorithm \ref{PropEitherOr}, we go to the point $$V_2=(\{3,6\}\{1,4\}\{2,5,7\}),$$
                                     which  is connected  with $V_1$ by an edge labeled by  $(\{3,6\}\{1,4\}\{7\}\{2,5\})$.
                                     Next come the vertices  $$V_3=(\{3,4,6\}\{1\}\{2,5,7\})$$
                                    $$  \hbox{  and    } V_4=(\{3,4,6\}\{1,
                                    2\}\{5,7\}).$$
                                  %  we cannot add more since $l_1+l_2+l_3>\frac{L}{2}=19\frac{1}{2}.$
                                      Then comes the vertex $$V_5=(\{3,4\}\{1, 2\}\{5,7,6\})=(\{4,3\}\{2,1\}\{7,6,5\}),$$
                                     which is the mirror image  of the target point.

                                      Now we start turning the configuration inside out, as is prescribed in Algorithm \ref{AlgINsideOut}.
                                     \item The next point is  $$V_6=(\{3,4\}\{1, 2,7,6\}\{5\}).$$
                                     After freezing the middle set, we get a triangular configuration of a quadrilateral
                                     to be turned inside out in three steps:
                                       $$V_7=(\{3\}\{1, 2,7,6\}\{5,4\}),$$
                                      $$V_8=(\{5,3\}\{1, 2,7,6\}\{4\}),$$
                                      $$V_9=(\{5\}\{1, 2,7,6\}\{3,4\})=(\{5\}\{6,7,1, 2\}\{3,4\}).$$
                                       One more edge brings us to the target point
                                      $$V_{10}=(\{5,6,7\}\{1, 2\}\{3,4\})=V'.$$
                                   \end{enumerate}
\end{ex}

\section{ Navigation and control on the moduli space}\label{SecNavigationOnModuli}
Here we describe a  finite-step algorithm of navigation from an
arbitrary point (which is not necessarily a vertex  of $\Gamma(L)$)
of the moduli space $M(L)$ to an arbitrary target  point.

We work under assumption that we can somehow control the shape of a
convex configuration. As an example, we can use the result of
\cite{ADEHOST}, where it is shown that a convex polygon can be moved
into any other convex polygon with the same counterclockwise
sequence of edge lengths in such a way that each angle changes
monotonically.

At the same time, we explain how to realize the flex explicitly.
As in the previous section, we assume that $l_1$ is the biggest
edgelength.
\bigskip

\bigskip

\textbf{Algorithm B}

\begin{enumerate}
 \item Given a starting configuration $S$ and a target configuration
 $T$, we take the $(n-3)$-cells of the complex $\mathcal{K}(L)$ containing $S$ and $T$.
 The cells may be not unique, if this is the case, choose $S$ and $T$ to be
 any of them.
 We choose $V_S=(I,\{1\},J),$ and $V_T$ to be some vertices of these
 two cells. Starting from now, we keep in mind Algorithm A applied
 for the vertices $V_S$ and $V_T$.

  \item
We  navigate  from $S$ to $V_S=(I,\{1\},J)$.
 In particular, this means that we spare one step compared to Algorithm \ref{PropEitherOr}.

We realize both $P$ and the convexification  $Conv(P)$   as two
bar-and-joint mechanisms.  By construction, there is a natural
bijection between edges of the two polygons, and the corresponding
edges are parallel. For each pair  of parallel edges (one edge from
$P$, and the other one from $Conv(P)$), we plug in a pair of
parallelograms  as is shown in Figure \ref{FigParall}.

\bigskip
\begin{figure}[h]
\centering
\includegraphics[width=8 cm]{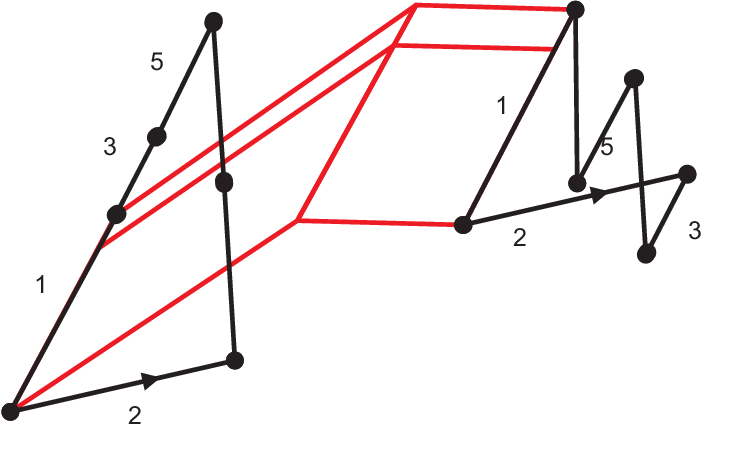}
\caption{Connecting parallelograms}\label{FigParall}
\end{figure}

To prevent turning inside out, we add one extra edge inside each of  the parallelograms (here we follow \cite{KM}).

Each  (convex) flex of $Conv(P)$ uniquely induces a flex of
 $P$ in such a way that the first polygon remains the convexification  of the second one. Therefore, the task is to bring the
 convex polygon $Conv(P)$ to a triangular configuration.
By assumption,  we can control $Conv(P)$, and therefore, we have a
controlled way of bringing $P$ to the vertex $(I,\{1\},J).$

  \item
On this step,  we  navigate  according to Algorithm A by prescribed
edges from one vertex $V$ to some other vertex $V'$.

We realize the motion  almost  in  the  same way as
above: Take any point $P$ lying on the  edge connecting $V$ and $V'$
and again realize both $P$ and the convexification   $Conv(P)$ as
bar-and-joint mechanisms.
 Now the polygon $Conv(P)$  is a quadrilateral. Each  of its edges we decompose into small edges of lengths $l_i$.
 Thus we again have a natural bijection between edges
of the two polygons, and the corresponding edges are parallel. For
each pair  of parallel edges (one edge from $P$, and the other one
from $Conv(P)$), we plug in  parallelograms in the same way as we
did above.

We can assume that the edges of the quadrilateral are frozen, that
is, the quadrilateral can flex only at the four vertices.

Each  (convex) flex of $Conv(P)$ uniquely induces a flex of
 $P$ in such a way that the first polygon remains the convexification  of the second one. Therefore, the task is to bring the
 convex polygon $Conv(P)$ from one triangular configuration to another  triangular configuration, see Figure \ref{FigEdge}.
 This can be realized in many ways, since   a quadrilateral has one degree of freedom, and its flexes  are well-understood, see
 \cite{KHPS}.

 Important is that every next edge  needs  a separate collection of auxiliary  parallelograms.
  \item Once we arrive at the point $V_T$, we go to the target point
  $T$ as on the very first step.
\end{enumerate}

\bigskip

In our second approach, we again add auxiliary bars to the polygonal
linkage, but now we  have one and the same bar-and-joint mechanism
which is not rearranged during the desired flex.

\begin{figure}[h]
\centering
\includegraphics[width=8 cm]{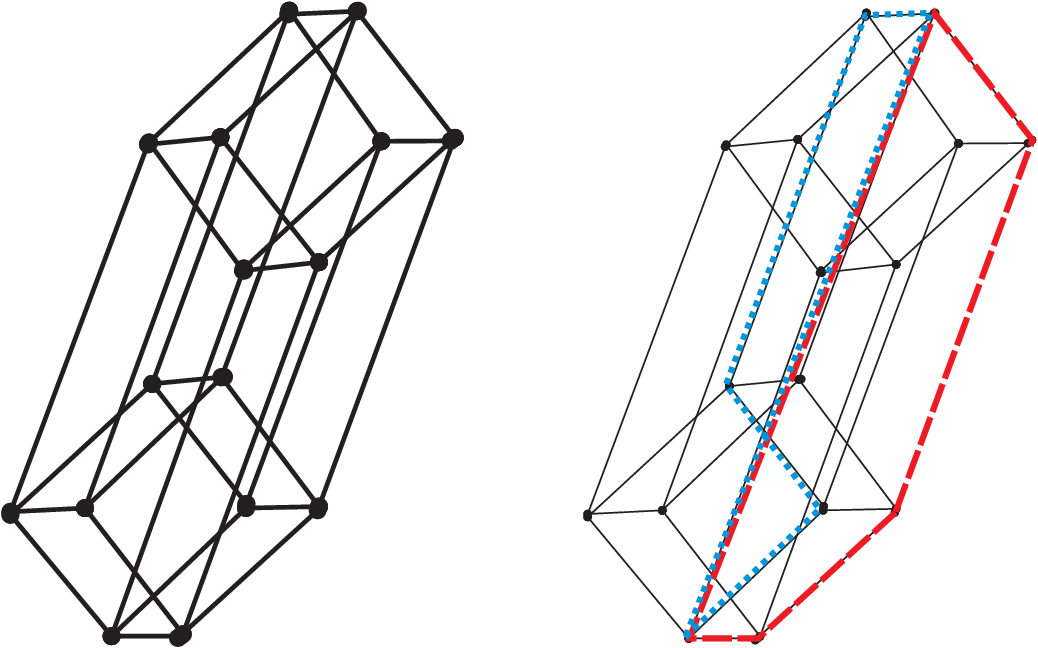}
\caption{Projection of the hypercube  includes both $P$ (blue) and
$Conv(P)$ (red). }\label{FigHypercube}
\end{figure}

The key observation is that the  projection of the 1-skeleton  of
hypercube can serve as a universal permuting machine: together with
any configuration $P$, it contains all other configurations obtained
from $P$ by permuting the order of edges, see Figure
\ref{FigHypercube}.

On the one hand, an obvious advantage of this approach is that we do
not remove and add bars at each step. On the other hand, a
disadvantage is that we add many extra bars.
Therefore, the choice between the two ways  depends on the
particular  task.

\bigskip

\textbf{Algorithm C}

\begin{enumerate}
  \item We assume that the starting and the target points are  $S, \ T\in M(L)$ . We find
  the vertices $V_S$ and $V_T$  exactly as in Algorithm B.
  \item
 Interpret the edges of $S$ numbered $1,2,...,(n-1)$  as projections of edges of the hypercube $[0,1]^{n-1}$. The edge
  numbered by $n$ is then the projection of the diagonal of the
  hypercube. Add extra bars to incorporate $S$ to the entire projection of the hypercube,
   which is now treated as a bar-and-joint mechanism, see Figure
\ref{FigHypercube}.
  \item The new mechanism includes also $Conv(S)$. Now we can manipulate by $Conv(P)$ following Algorithm A.
  As in Algorithm B, we first bring $Conv(S)$ to the  configuration labeled by $V_S=(I,\{1\},J).$
  \item Next, we follow the path on the graph $\Gamma$ prescribed by Algorithm A. For each step, we
  manipulate with the convex configuration $Conv(P)$ which degenerates to a convex quadrilateral.
  \item The last step  from $V_T$ to $T$ is the same as the very
  first step.
\end{enumerate}

\bigskip


\begin{thebibliography}{7}




\bibitem{ADEHOST} O. Aichholzer, E. D. Demaine,
J. Erickson, F. Hurtado, M. Overmars, M. Soss, G. Toussaint,
\textit{Reconfiguring convex polygons}, Computational Geometry 20
(2001), 85-95.

\bibitem{faS} M. Farber and D. Sch\"{u}tz, \textit{Homology of planar polygon
spaces,} Geom. Dedicata, 125 (2007),  75-92.

\bibitem{KM} M. Kapovich and J. Millson, \textit{Universality theorems for configuration spaces
of planar linkages,} Topology, 41 (2002), 1051-1107.

\bibitem{KapovichMillson} M. Kapovich and J. Millson, \textit{On the moduli space
of polygons in the Euclidean plane}, J. Diff. Geom., 42 (1995),
430-464.

\bibitem{LenWh} W. J. Lenhart, S. H. Whitesides, \textit{Reconfiguring closed polygonal
chains in Euclidean $d$-space,} Discrete and Computational Geometry
13 (1995),  Issue 1, 123-140.

\bibitem{LiuTr} G.F. Liu, J.C. Trinkle, \textit{Complete Path Planning for Planar Closed Chains Among Point
Obstacles}, in: Proceedings of Robotics: Science and Systems, 2005,  Cambridge, USA.

\bibitem{Pan} G. Panina,   \textit{ Moduli space of planar polygonal linkage: a combinatorial description}, arXiv:1209.3241


\bibitem{hausmann}J.-C. Hausmann, \textit{Controle des bras articules et transformations de Mobius.} L'Enseignement Mathematique, 51 (2005) 87-115.
\bibitem{hausmannrodriguez}J.-C. Hausmann and E. Rodriguez, \textit{Holonomy orbits of the snake charmer algorithm}, Geometry
and Topology of Manifolds, Banach Center Publications,  76  (2007)

\bibitem{khim2010}G. Khimshiashvili,  \textit{Maxwell problem for polygonal linkages}, Bull. Georgian Natl. Acad. Sci. 6, 2012, No. 2, 17-22.

\bibitem{khipan}G. Panina,  G. Khimshiashvili, \textit{On the area of a polygonal linkage},
Doklady Akademii Nauk,  Mathematics, 85, No. 1 (2012), 120-121.
\bibitem{KHPS}{G. Khimshiashvili, G. Panina, and  D. Siersma, }\textit{Coulomb control of polygonal
linkages,} J. of Dynamical and Control Systems, 20, 2014, No. 4,
491-501.

\end{thebibliography}
\end{document}